\documentclass[12pt]{amsart}
\usepackage{amsfonts,amsthm}

\newtheorem{theorem}{Theorem}[section]
\newtheorem{corollary}[theorem]{Corollary}
\newtheorem{lemma}[theorem]{Lemma}
\newtheorem{proposition}[theorem]{Proposition}
\newtheorem*{acknow}{Acknowledgments}
\newtheorem{remark}[theorem]{Remark}
\theoremstyle{definition}
\newtheorem*{defn}{Definition}
\theoremstyle{remark}

\newtheorem{case}{Case}

\newcommand{\nc}{\newcommand}
\newcommand{\be}{\begin{equation}}
\newcommand{\ee}{\end{equation}}
\newcommand{\bea}{\begin{eqnarray}}
\newcommand{\eea}{\end{eqnarray}}
\newcommand{\no}{\nonumber}

\newcommand{\bc}{\begin{center}}
\newcommand{\ec}{\end{center}}

\nc{\benu}{\begin{enumerate}}
\nc{\eenu}{\end{enumerate}}

\nc{\bth}{\begin{theorem}}
\nc{\eth}{\end{theorem}}
\nc{\bpr}{\begin{proposition}}
\nc{\epr}{\end{proposition}}
\nc{\ble}{\begin{lemma}}
\nc{\ele}{\end{lemma}}
\nc{\bco}{\begin{corollary}}
\nc{\eco}{\end{corollary}}
\nc{\bre}{\begin{remark}}
\nc{\ere}{\end{remark}}
\newcommand{\la}{\lambda}

\nc{\f}{\frac}
\nc{\rw}{\rightarrow}
\nc{\Rw}{\Rightarrow}
\newcommand{\st}[1]{\hat #1}
\nc{\De}{\Delta}

\nc{\na}{\nabla}
\nc{\al}{\alpha}
\nc{\rai}{R(e_\al ,e_i ,e_\al ,e_i)}
\nc{\bet}{\beta}
\nc{\rab}{R_{\al \bet \al \be}}
\nc{\ga}{\gamma}
\nc{\de}{\delta}
\nc{\ep}{\varepsilon}
\nc{\ei}{\ep _i}
\nc{\ea}{\ep _{\al}}
\nc{\ej}{\ep _j}
\nc{\eb}{\ep _{\bet}}
\nc{\br}{\bgroup\em}
\nc{\er}{\egroup}
\nc{\rs}{Riemannian submersion}
\nc{\srs}{semi-Riemannian submersion}
\nc{\tuf}{totally umbilic fibres}
\nc{\co}{compact and orientable total space}
\nc{\com}{compact and orientable manifold}
\nc{\om}{\omega}
\nc{\ta}{\tau ^{HV}}

\newcommand{\pis}{$\pi :(M,g) \to (B,g')$}
\nc{\hauu}{h^\al _{11}}
\nc{\haij}{h^\al _{ij}}
\nc{\aiab}{A^i _{\al \bet}}
\nc{\dvg}{\mathrm{dv}_g}
\nc{\ve}{$\mathcal V$}
\nc{\ho}{$\mathcal H$}
\nc{\tm}{\Gamma (TM)}
\nc{\tb}{\Gamma (TB)}
\nc{\lc}{Levi-Civita connection}
\nc{\dec}{\st{\de}}
\nc{\deu}{\check\de}
\begin{document}
\title{Semi-Riemannian submersions with totally umbilic fibres}
\author{Gabriel B\u adi\c toiu}
\author{Stere Ianu\c s}
\address{Institute of Mathematics             
    of the Romanian Academy,              
    P.O. Box 1-764,                      
    Bucharest 70700,                      
    Romania}
\email{gbadit@stoilow.imar.ro}
\address{University of Bucharest,
    Department of Mathematics, C.P. 10-119, Post. Of. 10,
    Bucharest 72200,
    Romania}
\email{ianus@pompeiu.imar.ro}
\date{}
\maketitle
\section*{Introduction}

The theory of Riemannian submersions was initiated
by O'Neill \cite{one} and Gray \cite{gra}.
Presently, there is an extensive literature on the \rs s 
with different conditions imposed on the total
 space and on the fibres.
A systematic exposition could be found in the A.Besse's book \cite{bes} .
Semi-\rs s were introduced by O'Neill in his book \cite{onei}; for Lorentzian case see Magid
\cite{mag}.
In this paper we study the \srs s with \tuf .\\
The main purpose of section \S 2 is to obtain obstructions to
the existence of the \srs s  with \tuf\ and with compact and orientable total
space, in terms of sectional and scalar curvature.
Using a formula of Ranjan \cite{ranj}, we obtain an integral 
formula for mixed scalar curvature $\ta$,
 which give us obstructions to existence
of the \srs s in some special cases.
Then we establish other integral formula of scalar curvature
of total, base and fibres spaces and
another obstruction to existence of \srs s  is obtained.\\
In section \S 3 we study the semi-Riemannian submersions
$\pi:M\to B$ with totally umbilic fibres,
when  the mean curvature vector field $H$
 is parallel in the horizontal bundle along fibres and
 $R(X,Y,X,Y)$ is constant along fibres
 for every $X$, $Y$ basic vector fields.
If moreover we assume that
$B$ is a Riemannian manifold
and $M$ is a semi-Riemannian manifold of index 
$\mathrm{dim\ }M-\mathrm{dim\ }B$
we deduce that the mean curvature vector field $H$ is basic if and only if
the horizontal projection of $R(X,Y)A_XY$, denoted by 
$hR(X,Y)A_XY$ (see page 2), is basic  for every $X$, $Y$ basic vector fields.
We get there are no semi-Riemannian submersions $\pi:M\to B$ 
with totally umbilic fibres, $M$ 
a constant positive curvature semi-Riemannian manifold of index 
$\mathrm{dim\ }M-\mathrm{dim\ }B\geq 2$
and $B$ a compact and orientable
Riemannian manifold.
Then we find that a semi-Riemannian submersion 
with totally umbilic fibres,
with $R(X,Y,X,Y)$ constant along fibres for every $X$, $Y$ basic vector fields
and 
with basic mean curvature vector field
from an $m$-dimensional semi-Riemannian manifold
of index $r=m-n$ with non-negative mixed curvature
onto an $n$-dimensional compact and orientable Riemannian manifold,
has totally geodesic fibres, integrable horizontal distribution
and null mixed curvature.
Therefore a semi-Riemannian submersion
$\pi:M\to B$ with totally umbilic fibres and 
with sectional curvature of the fibres non-vanishing anywhere,
from a constant curvature
semi-Riemannian manifold $M$ of index 
$\mathrm{dim\ }M-\mathrm{dim\ }B\geq 2$  is a Clairaut
semi-Riemannian submersion.
Also we study the case of positive sectional  curvature  fibres.
We give a sufficient condition to have every fibre with 
zero sectional curvature, when the total space has
constant curvature.

\section{Preliminaries}

In this section we recall some notions and results which will be needed.
\begin{defn}
    Let $(M,g)$ be an $m$-dimensional connected 
    semi-Riemannian manifold of index $s$ ($0 \le s\le m$),
    let $(B,g')$  be an $n$-dimensional connected 
    semi-Riemannian manifold of index $s'\le s$, ($0 \le s'\le n$).
    {\it A semi-Riemannian submersion} (see O'Neill \cite{onei}) is 
    a smooth map
         $\pi :M\to B$
    which is onto and satisfies the following three axioms:
\begin{itemize}
\item[(a)] $\pi_*|_p$ is onto for all $p\in M$;
\item[(b)] The fibres $\pi^{-1}(b)$,  $b\in B$ are 
    semi-Riemannian submanifolds of $M$;
\item[(c)] $\pi_*$ preserves scalar products of 
    vectors normal to fibres.
\end{itemize}
\end{defn}

  We shall always assume that the dimension of the fibres
  $\dim M-\dim B$ is positive and
  the fibres are connected.

    The vectors tangent to fibres are called 
    vertical and those normal to fibres
    are called horizontal.
    We denote by \ve\ the vertical distribution and by \ho\
    the horizontal distribution.\\
    B.O'Neill \cite{one} has characterized the geometry 
    of a Riemannian submersion in terms of
    the field tensors $T$ , $A$ defined for 
    $E$ , $F \in \tm $ \ by
\bea   A_E F&=&h \na _{hE}{vF} + v \na _{hE}{hF}\no\\
       T_E F&=&h \na _{vE}{vF} + v \na _{vE}{hF}\no
\eea
    where $\na $ is the \lc\ of $g$.
    Here the symbols $v$ and $h$ are the orthogonal 
    projections on \ve\ and \ho\
    respectively. For basic properties of \rs s and examples see
    \cite{bes}, \cite{gra}, \cite{one}.
    The letters $U$, $V$, $W$, $W'$ will always denote 
    vertical vector fields,
     $X$, $Y$, $Z$, $Z'$ horizontal vector fields 
    and $E$,$F$,$G$,$G'$
    arbitrary vector fields on $M$.
    A vector field $X$ on $M$ is said to be basic if 
    $X$ is horizontal and $\pi -$related to a vector 
    field $X'$ on B.
    It is easy to see that every vector field $X'$
    on $B$ has a unique horizontal lift $X$ to $M$ 
    and $X$ is basic. 
    The following lemmas are well known
    (see \cite{one},\cite{esc}). 
\ble\label{lem}
      Let $X$ be a horizontal vector field. 
    If $g_p(X,Z)=g_{p'}(X,Z)$
    for all $Z$ basic vector fields on $M$, 
    for all $p$,$p' \in \pi ^{-1}(b)$
    and for all $b \in B$ then $\pi _*X$ is a well 
    defined vector field on
    $B$ and $X$ is basic.
\ele
\ble\label{lemdoi}
    We suppose X and Y are basic vector fields on M
    which are $\pi$-related to $X'$ and
    $Y'$, and $V$ is a vertical vector field.
    Then \\
    a) h\br [\er X,Y \br ]\er is basic and
     $\pi$-related to $[X', Y']$;\\
    b)\ $h\na _X Y$ is basic and $\pi$-related to 
    $\na '_{X'}Y'$\ , where $\na '$ is the \lc\ on B;\\
    c) $h\na_VX=A_XV.$
\ele
   Let $\{ e_1$,...,$e_m\}$ be a local field of orthonormal
   frames on $M$ such that $e_1,...,e_r$ are vertical vector fields
   and $e_{r+1},...,e_m$ are basic vector fields, 
   where $r=\mathrm{dim\ }M-\mathrm{dim\ }B$ denotes the dimension of fibres.
    We have $g(e_a,e_b)=\ep _a\de _{ab}$ for every $a$, $b$ 
    (where $\ep _a\in\{-1,1\}$).
    We shall always denote {\it vertical} indices by $i,j,k,l,... =1,...,r$
    and {\it horizontal} indices by $\al , \bet , \ga , \de ,... =r+1,...,m$.
    The summation $\sum $\ is taken
    over all repeated indices, unless otherwise stated.\\
    The convention for the Riemann tensor used is 
     $R(E,F)G=\na_E\na_FG-\na_F\na_EG-\na_{[E,F]}G$ and 
     $R(E,F,G,G')=-g(R(E,F)G,G')$.\\
    Let $\st{g}$ be the semi-Riemannian metric
    of a fibre $\pi^{-1}(b)\ ,\ b\in B$.
    We make the following notations:\\
       $R,R',\st{R}$\ for the Riemann tensors,
       $K,K',\st{K}$\ for the sectional curvatures, and
       $s,s',\st{s}$\ for the  scalar
    curvatures of the metrics $g,g',\st{g}$, respectively;
\bea
     \ta &=& \sum \ei\ea\rai\ ,\ H\ =\ \sum \ei T_{e_i}e_i,\no\\ 
     g(A,A) &=& \sum\limits_{\al,\bet}\ea\eb g(A_{e_\al}e_\bet,A_{e_\al}e_\bet)
            =\sum\limits_{\al , i}\ei\ea g(A_{e_\al}e_i,A_{e_\al}e_i)\ ,\no\\
     g(T,T) &=& \sum\limits_{i , j}\ei\ej g(T_{e_i}e_j,T_{e_i}e_j)
            =\sum\limits_{\al,i}\ei\ea g(T_{e_i}e_\al,T_{e_i}e_\al)\ ,\no\\
     \mathrm{div}(E) &=& \sum\limits_i \ei g(\na _{e_i}E,e_i)+
                \sum\limits_\al \ea g(\na _{e_\al}E,e_\al).\no
\eea
If $X$ is  an unitary horizontal vector and
$V$ is an unitary  vertical
vector the sectional curvature of 
the $2$-plane $\{X,V\}$ is called
the {\it mixed sectional curvature}.
\begin{defn}
    A semi-Riemannian submanifold $F$ of a semi-Riemannian manifold
    $(M,g)$ is said to be {\it totally umbilic
    submanifold} if the second fundamental form $\Pi$ of $F$ is given by
    $\Pi(U,V)=g(U,V)\f Hr$ for every $U$, $V$ tangent vector fields to $F$.
\end{defn}

    Notice that $T_UV$ is the second fundamental form of the fibres
    and $A_XY$ is a natural {\it obstruction} to integrability of horizontal distribution.

    If the fibres of the semi-Riemannian submersion are totally
    umbilic submanifolds then $T_UV=\f {1}{r}g(U,V)H$ for every $U$, $V$ vertical
    vectors fields and $g(T,T)=\f {1}{r}g(H,H)$.

    By O'Neill's equations \cite{one} we get the following lemma.
\ble\label{lemtrei}
    If \pis\ is a \srs\ with totally umbilic fibres then:
\bea
    a)\ R(U,V,U,V)&=&\st{R}(U,V,U,V) +[g(U,V)^2-g(U,U)g(V,V)]g(\f Hr,\f Hr)\no\\
    b)\ R(X,U,X,U)&=&g(U,U)[g(\na_X\f Hr,X)-g(X,\f Hr)^2]+g(A_XU,A_XU)\no\\
    c)\ R(X,Y,X,Y)&=&R'(\pi_*X,\pi_*Y,\pi_*X,\pi_*Y)-3g(A_XY,A_XY)\no
\eea
\ele
    Using a relation of R. Escobales and Ph. Parker \cite{escob}, 
    we have the following proposition
\bpr\label{unupatru}
    Let \pis\ be a 
    semi-Riemannian submersion with totally umbilic fibres and
    $X$, $Y$ be basic vector fields.
    Then $A_XY$ is a Killing vector field along fibres if and only if 
    $g(\na_YH,X)=g(\na_XH,Y)$.
\epr
    \begin{proof}
    Let $U$,$V$ be vertical vector fields.
    Since the fibres are totally umbilic, for every 
    $X$, $Y$ basic vector fields we have (see \cite{escob})
    $$g(\na_U(A_XY),V)+g(\na_V(A_XY),U)=
  	  \f{g(U,V)}r (g(\na_YH,X)-g(\na_XH,Y)).$$
\end{proof}
\bco\label{unucinci}
    Let \pis\ be a 
    semi-Riemannian submersion with totally umbilic fibres.
    We suppose that the mean curvature 
    $H$ is a basic vector field. 
    Then $A_XY$ is a Killing vector field along fibres for every $X$,$Y$ basic
    vector fields if and only if $(\pi_*H)^\flat$ is a closed $1-$form on
    $B$.
\eco 
\begin{proof}
    Let $\om =(\pi_*H)^\flat$, $X'$,$Y'$ be vector fields on $B$.
    We have
     $$2\mathrm{d}\om)(X',Y')=g'(\na_{X'}\pi_*H,Y')-g'(\na_{Y'}\pi_*H,X').$$
    Let $X$,$Y$ be basic vector fields such that $\pi_*X=X'$, $\pi_*Y=Y'$.
    By lemma \ref{lemdoi}, we have 
    $$(g'(\na_{X'}\pi_*H,Y')-g'(\na_{Y'}\pi_*H,X'))\circ\pi
                 =g(\na_XH,Y)-g(\na_YH,X).$$
    Applying proposition \ref{unupatru} we get the conclusion.
\end{proof}
    In \cite{bad} we proved the following result.
\bpr\label{simply}
    Let \pis\ be a semi-Riemannian submersion 
    from an $(n+s)$-dimensional semi-Riemannian
    manifold of index $s\geq 1$ onto an $n$-dimensional Riemannian manifold.
    If $M$ is geodesically complete and simply connected then
\begin{itemize}
\item[$1)$] $B$ is complete and simply connected;
\item[$2)$] we have an exact homotopy sequence 
       $$\cdot\cdot\cdot\to\pi_2(B)\to\pi_1(fibre)\to\pi_1(M)\to\pi_1(B)\to 0;$$
\item[$3)$] If moreover $B$ has non-positive curvature then the fibres are simply 
       connected.
\end{itemize} 
\epr
    The following proposition is a semi-Riemannian version of Ranjan's formula 
(see \cite{ranj}) in the case with totally umbilic fibres.
\bth\label{forjan}
      If the \srs\ \pis\ has totally umbilic fibres, then
      \begin{equation}\label{form}
          \ta =\mathrm{div}(H)+(1-\f 1r)g(H,H)+g(A,A)
      \end{equation}
\eth
\begin{proof}
    From O'Neill \cite{one}, we have the following formula
    $$  \rai =g(\ (\na _{e_\al }T)_{e_i}e_i,e_{\al })
      +g(A_{e_\al }e_i,A_{e_\al }e_i)-
      g(T_{e_i}e_\al ,T_{e_i}e_\al )
   $$
  $
    g(\ (\na _{e_\al }T)_{e_i}e_i,e_{\al })
   =g(\na _{e_\al }(T_{e_i}e_i),e_{\al })
   -2g(T_{e_i}(v\na  _{e_\al }e_i),e_{\al })
  $
  
    If we denote by
    $q^{\al}_{ij}=g(v\na  _{e_\al }e_i,e_j)$, 
  it is easy to see that $q^{\al}_{ji}+q^{\al}_{ij}=0 $.
  Since $T_{e_i} e_j=T_{e_j} e_i$\ we get
  $\sum\limits_{j} \ej g(T_{e_j}(v\na  _{e_\al }e_j),e_{\al })=
  \sum\limits_{i , j}\ei \ej q^{\al}_{ji}T_{e_j}e_i =0. $\\
  Then
  \bea
     \ta &=&\sum\limits_{\al , i}\ei \ea \rai =\sum\limits_{\al }
     \ea g(\na _{e_\al }(\sum\limits_i\ei T_{e_i}e_i),e_{\al })+\no\\
     &&\sum\limits_{\al , i}\ea \ei g(A_{e_\al }e_i,A_{e_\al }e_i)
     -\sum\limits_{\al , i}\ea \ei g(T_{e_i}e_\al ,T_{e_i}e_\al )\no
  \eea
  We get the semi-Riemannian version of
Ranjan's formula
     $$\ta =\mathrm{div}(H)+g(H,H)+g(A,A)-g(T,T).$$
  Since the  fibres of the \srs\  $\pi$ are totally umbilic
  we have\\ 
        $g(T,T)=\f 1rg(H,H)$.
\end{proof}

\section{Integral formulae}
As a consequence of theorem \ref{forjan}, we have the following
proposition.
\bpr\label{for}
   If the \srs\ \pis\ has totally umbilic fibres
   and M is a \com , then 
      \begin{equation}\label{intfor}
        \int\limits_M \ta \dvg
           =(1-\f 1r)\int\limits_{M}g(H,H)\dvg +
         \int\limits_{M}g(A,A)\dvg
     \end{equation}
\epr
\begin{proof} We use the relation (\ref{form}) and $\int\limits_M
     \mathrm{div}(H)=0$. 
\end{proof}
\bco\label{thdoi}
     If $\pi :M\to B$\ is a \rs\ with \tuf\ and M is
      a \com\ then $\int\limits_M\ta \dvg \ge 0$.
\eco
     By proposition \ref{for}, we have the following
     splitting theorem, which is a generalization
     of proposition 3.1. of R. Escobales \cite{esc}.
\bth\label{t:2.3}
       Let $\pi :M\to B$ be a Riemannian submersion
       with \tuf . We suppose that M is a \com\ with 
       non-positive mixed sectional curvature 
     $($ i.e. $K(X,V)\leq 0$ for every $X$ horizontal vector
    field and for every $V$ vertical vector field $)$. Then 
\begin{itemize}
\item[a)] $R(X,U,Y,V)=0$ for every $X$, $Y$ horizontal vector fields and for
       every $U$ , $V$ vertical vector fields;      
\item[b)] the horizontal distribution is integrable
         (this is equivalent with $A \equiv 0$);
\item[c)] the fibres are totally geodesic.
\end{itemize}
\eth
\begin{proof} 
      Since $\ta \le 0$\ ,\ $g(A,A)\geq 0$ and $g(H,H)\geq 0$,
      we get the relations
       $\ta=\sum\limits_{i\ ,\al}\rai= 0$
       \ ,\ $g(A,A)=0$ and $(1-\f 1r)g(H,H)=0$, by
      formula (\ref{intfor}).
      Hence $A \equiv 0$.\\
      Let $X$ be a horizontal vector and
      $V$ be a vertical vector, $X\neq 0$, $V\neq 0$.
      We can choose  a local field $e_1,...,e_m$ of orthonormal
      frames adapted to the Riemannian submersion such that
      $e_1=\f V{\| V\|}$, $e_{r+1}=\f X{\| X\|}$.
      Since $\ta =0$ and $\rai\leq 0$ for all $i, \al$ we get $\rai =0$ 
       for all $i,\ \al$.
      Therefore $R(X,V,X,V)=0$ for every
      $X$ horizontal vector field and for every $V$ vertical vector field.

     We shall prove that the fibres
     are totally geodesic.\\
     Since $A\equiv 0$, $R(X,V,X,V)=0$, we have, by lemma \ref{lemtrei},
\begin{equation}\label{e:e}
     g(\na_X\f Hr,X)=g(X,\f Hr)^2
\end{equation}
     for every $X$ horizontal vector field.
     Let $p$ be an arbitrary point in $M$ and $\ga:\mathbb R\to M$
     a horizontal geodesic in $M$, $\ga(0)=p$. 
     We denote by $h(t)=g(\f Hr,\dot\ga(t))$. Rewriting formula \eqref{e:e},
     for every $t\in\mathbb R$ we get
\begin{equation}\label{e:dif}
     \f{\mathrm dh}{\mathrm dt}(t)=h(t)^2
\end{equation}
     The differential equation \eqref{e:dif} has the solution
     $h(t)\equiv 0$ or $h(t)=-\f 1{t+A}$ for some $A\in\mathbb R$.
     But the domain of the maximal solution is the entire real line,
     only for the null solution. Hence $g(\f Hr,\dot\ga(t))\equiv 0$
     for every $\ga$ horizontal geodesic. Therefore
     $H\equiv 0$.
\end{proof} 
\bpr\label{p:2.4}
     Let $B$ be an $n$-dimensional Riemannian manifold,
     $\pi :M\to B$ be a \srs\ with \tuf\ and M be an 
     $m$-dimensional compact, orientable semi-Riemannian 
     manifold of index $r=m-n$.
     Then
     \begin{equation}\label{ftreia}
        \int\limits_M(  \ta -(1-\f 1r)g(H,H)  )\dvg \le 0,
     \end{equation}
     we have equality in \eqref{ftreia} if and only if the horizontal distribution is
     integrable;
     \begin{equation}\label{ftreib}
        \int\limits_M(  \ta -g(A,A)  )\dvg \geq 0,
     \end{equation}
     we have equality in \eqref{ftreib} if and only if either $r=1$ or the fibres are
     totally geodesic.
\epr
\begin{proof} Since $\pi $ sends isometrically the horizontal spaces
        into the tangent space of $B$, 
        and $B$ is a Riemannian manifold,
        it follows that the fibres are semi-Riemannian manifolds of
        indices $r$ and $g(H,H)\geq 0$.
        Since $g(A,A)=\sum\limits_{\al , \bet}\ea\eb 
        g(A_{e_\al}e_\bet ,A_{e_\al}e_\bet)$, $\ea =1$ for all
        $\al$ , and
        the induced metrics on fibres are negative definite, 
        we obtain $g(A,A)\le 0$. 
        By formula (\ref{intfor}) we have the conclusion. 
\end{proof}
Applying proposition \ref{p:2.4} for a Lorentzian total space we get the following
obstruction
.\bth
     Let $B$ be a $n$-dimensional Riemannian manifold. 
     If $M$\ is a compact, orientable $(n+1)$-dimensional 
     semi-Riemannian space of index $1$, with positive mixed curvature
     then there are no $\pi :M\to B$ \srs s.
\eth
\begin{proof} 
       We suppose that there is such a semi-Riemannian submersion.
       Since $r=1$ the inequality (\ref{ftreia}) implies
       $\int\limits_M \ta \le 0$, which is
       a contradiction with the stated condition of positive mixed curvature.
\end{proof}
\bpr
      If $M$ is a $(r+1)$-dimensional semi-Riemannian manifold
      of index $r$, with negative mixed curvature and 
      $B$ is a one dimensional Riemannian
      manifold then there are no $\pi:M\to B$ \srs s 
      with compact, orientable total space and \tuf .
\epr
\begin{proof}
       We suppose that there is such a semi-Riemannian submersion.
       Since dim \ho\ =1 we have $A\equiv 0$ .
       By relation (\ref{ftreib}) we obtain $\int\limits_M \ta\geq 0$,
       which is a contradiction with negative mixed curvature condition.
\end{proof}
\begin{defn}
   The {\it total scalar curvature} of a compact manifold
   $(M,g)$ is $S_g=\int_Ms_g\dvg$, 
   where $s_g$ is the scalar curvature
   of $(M,g)$.
\end{defn}

  In what follows we give some obstructions to the
  existence of \srs s in terms of
  total scalar curvatures of the total space $M$, 
  the base space $B$ and of the fibres.

  We denote by $S$, $S'$, $\hat S$ the total scalar curvature
   of $(M,g)$, $(B,g')$ and $(\pi^{-1}(x),\hat g)$ respectively.
   
  By lemma \ref{lemtrei}  and formula \eqref{form},
  we get immediately the following proposition.
\bpr
     If  $\pi:M\to B$ is a \srs\ with \tuf ,
      then 
    \begin{equation}
        \label{scal}
         s=s'\circ\pi +\st{s} +2\mathrm{div}(H) +(1-\f 1r)g(H,H)-g(A,A)
    \end{equation}
\epr
Integrating formula \eqref{scal} we get
\bpr
     If $\pi:M\to B$ is a \srs\ with \co\ M and with \tuf , then 
       \begin{equation}\label{sca} 
             S-(S'\circ\pi +\st{S})=(1-\f 1r)
            \int\limits_{M}g(H,H)\dvg -
            \int\limits_{M}g(A,A)\dvg
       \end{equation}
\epr
\bco
       Let $B$ be an $n$-dimensional Riemannian manifold,
       $\pi :M\to B$ be a \srs\ with \tuf\ and M be an 
       $m$-dimensional compact and
       orientable semi-Riemannian manifold of index $r=m-n$.
       Then 
           $S\ge S'\circ\pi +\st{S}.$ 
       We have equality if and only if the horizontal distribution is 
       integrable and either $r=1$ or
       the fibres are totally geodesic.
\eco
\begin{proof} 
            We have $g(H,H)\geq 0$,  $g(A,A)\le 0$. Hence, by formula
            (\ref{sca}), we get the conclusion.
\end{proof}
\bpr
      If $\pi:M\to B$ is a \rs\ with \co\ M and with fibres of
       dimension $1$ then
          $S\leq S'\circ\pi$
        We have equality if and only if the horizontal 
        distribution is integrable.
\epr
\begin{proof} We have $g(A,A)\ge 0$,
              $S-(S'\circ\pi +\st{S})=
              -\int\limits_M g(A,A)\dvg $ and $\st{s}=0$.
             Hence, by \eqref{sca}, we get $S\leq S'\circ\pi$.
\end{proof}

    Let $s^{\mathcal H}=\sum\limits_\al\ea\rho (e_\al,e_\al)$, where $\rho$ is the
  Ricci tensor of $M$ (see \cite{kwo}).
\bpr
     If \pis\ is a semi-Riemannian submersion with totally umbilic fibres then
\begin{equation}\label{rh}
     s^{\mathcal H}-s'\circ\pi=\mathrm{div}(H)+(1-\f 1r)g(H,H)-2g(A,A)
\end{equation}
\epr
\begin{proof}
     By lemma \ref{lemtrei} we have
\bea
    s^{\mathcal H}-s'\circ\pi 
      &=& 
        \sum\limits_{\al ,i}\ea\ei R(e_\al,e_i,e_\al,e_i)
        +\sum\limits_{\al ,\bet}\ea\eb R(e_\al,e_\bet,e_\al,e_\bet)\no\\
      &=& \sum\limits_\al\ea (g(\na_{e_\al}H,e_\al)-\f 1rg(e_\al,H)^2)+g(A,A)-
         3g(A,A)\no\\
      &=&  \mathrm{div}(H)-\sum\limits_i\ei g(\na_{e_i}H,e_i)-\f 1r
         g(H,\sum\limits_\al\ea g(e_\al ,H)e_\al)-2g(A,A)\no\\
      &=& \mathrm{div}(H)+(1-\f 1r)g(H,H)-2g(A,A)\no
\eea
\end{proof}

Integrating formula \eqref{rh} we get
\bpr
     If \pis\ is a semi-Riemannian submersion with totally umbilic fibres 
     and if $M$ is a compact and orientable manifold then
\begin{equation}\label{e:int}     
    \int\limits_M (s^{\mathcal H}-s'\circ\pi)\dvg 
         =(1-\f 1r)\int\limits_M g(H,H)\dvg -2\int\limits_M g(A,A)\dvg 
\end{equation}
\epr
\bco
     Let $B$ be an $n$-dimensional Riemannian manifold,
     $\pi :M\to B$ be a semi-Riemannian submersion 
     with totally umbilic fibres and M be an
     $m$-dimensional compact and orientable semi-Riemannian 
     manifold of index $r=m-n$.
     Then
      $$\int\limits_M (s^{\mathcal H}-s'\circ\pi)\dvg\geq 0$$
     We have equality if and only if the horizontal distribution is integrable
     and either $r=1$ or the fibres are totally geodesic.
\eco
\begin{proof}
    We have $g(H,H)\geq 0$ and $g(A,A)\leq 0$. Therefore, by \eqref{e:int}, we get the conclusion.
\end{proof}
\section{Mean curvature vector field}
  We denote by $\rho$, $\rho '$ and $\st{\rho}$ the Ricci tensors of the
  manifolds $M$, $B$ and of the fibre $\pi ^{-1}(b)$, $b \in B$.
  The letters $U$, $V$ denote vertical vector fields and $X$, $Y$
  horizontal vector fields.\\
  We introduce the following notations (cf. \cite{bes}):   
     \bea
        g(A_X,T_U)&=&\sum\limits_\al\ea g(A_Xe_\al,T_Ue_\al)=
                   \sum\limits_i\ei g(A_Xe_i,T_Ue_i),\no\\
        \deu A&=&-\sum\limits_\al\ea {(\na _{e_\al}A)}_{e_\al} \ , \ 
        \dec T=-\sum\limits_i\ei {(\na _{e_i}T)}_{e_i}\no
     \eea
    where $\ei\de _{ik} =g(e_i,e_k)$ for all $i , k $ , 
      $\ea\de _{\al\bet} =g(e_\al,e_\bet)$ for all $\al, \bet$. 
\bpr\label{tru}
    Let \pis\ be a \srs\ with totally umbilic fibres.
    Assume that the dimension of fibres $r\geq 2$. Then
    the following conditions are equivalent:
\begin{itemize}
\item[    i)\ ] the mean curvature vector field $H$ is a basic vector field;
\item[    ii)\ ] $g((\dec T)U,X)-g(A_X,T_U)=0$ for every $U$ vertical vector field
         and for every $X$ horizontal vector field;
\item[    iii)\ ] $\rho(X,U)+g((\deu A)(X),U)+(r+1)g(T_U,A_X)=0$ for every $U$
         vertical vector field
         and for every $X$ horizontal vector field.
\end{itemize}
\epr
\begin{proof} 
i)$\Leftrightarrow$\ ii)\\
    Since the fibres are totally umbilic we have 
     $g((\dec T)U,X)=-\f 1rg(\na_UH,X)$ 
     and $g(A_X,T_U)=\f 1rg(A_XU,H)$ for every $U$
     vertical vector field
     and for every $X$ horizontal vector field.\\
     Let $X$ be a basic vector field. By lemma \ref{lemdoi}, we have $A_XU=h\na_UX$.
     Then\\
     $g((\dec T)U,X)-g(A_X,T_U)=-\f 1rg(\na_UH,X)-\f 1rg(h\na_UX,H)=
                   -\f 1rUg(H,X)$.\\
    Hence, by lemma \ref{lem},
    $H$ is a basic vector field if and only if
    the function $g(H,X)$ is constant along fibres 
    for every $X$ basic vector field.\\
i)$\Leftrightarrow$\ iii)\\
    Let $X$ be a basic vector field.
    Using (9.36b) in \cite{bes}, we compute
  \bea   
    \rho(U,X) &=& g((\dec T)U,X)+g(\na_UH,X)-g((\deu A)(X),U)-2g(T_U,A_X) \no\\
              &=&
    g((\dec T)U,X)-g(T_U,A_X)+Ug(H,X)-g(H,A_XU)\no\\
     &&  
     -g((\deu A)(X),U)-\f 1r g(H,A_XU) \no\\
              &=&
    (1-\f 1r)Ug(H,X)-\f {r+1}r g(H,A_XU)-g((\deu A)(X),U).\no
  \eea
Then
   \begin{equation}\label{treizece}
           \rho(X,U)+g((\deu A)(X),U)+(r+1)g(T_U,A_X)=(1-\f 1r)Ug(H,X)
   \end{equation}
  Applying again lemma \ref{lem}, we get the conclusion.
\end{proof}
\bre
    By proposition \ref{tru} and theorem 2.2 in \cite{badi}, 
    $H$ is a basic vector field if and only if the contractions
    $(1,3)$ and $(2,4)$ of the Riemann tensor of the Vr\^ anceanu connection
    are equal.
\ere 

The following proposition is well known 
(see Proposition 9.104 in \cite{bes}).
\bpr\label{rem}
    If $B$, $F$ are semi-Riemannian manifolds
    and $f:B\to \Bbb R$ is a positive
    function then $\pi: B\times _fF\to B$ is a
    semi-Riemannian submersion with
    totally umbilic fibres, $A\equiv 0$ and $H$ a basic vector field.\\
   Conversely, if \pis\ is a semi-Riemannian  
   submersion with totally umbilic fibres 
   $A\equiv 0$ and $H$ a basic vector field then $M$ is a locally warped
   product.
\epr   
\bco\label{warped}
   If \pis\ is a semi-Riemannian submersion with totally umbilic fibres,
   if $M$ is an Einstein manifold, $r\geq 2$ 
   and the horizontal distribution \ho\ is integrable then
   $M$ is a locally warped product.
\eco
\begin{proof}
   Since $A\equiv 0$, $\rho(X,U)=0$, we have
   $H$ is a basic vector field,  by proposition \ref{tru}.
   By proposition 9.104 in \cite{bes} we have 
   $M$ is a locally warped product. 
\end{proof}

\begin{defn}
    A semi-Riemannian submanifold is said to be an {\it extrinsic sphere} if it is
    totally umbilic and the mean curvature vector field 
    $H$ is nonzero anywhere and
    parallel in the normal bundle.
\end{defn}

\bco
    Let $\pi:(M,g)\to (B,g')$
    be a semi-Riemannian submersion.
    If the mean curvature vector field $H$ is basic, 
    $\mathrm{dim\ }B=2$, and the fibres are extrinsic spheres
    then the horizontal distribution \ho\ is
    integrable and $M$ is a locally warped product.
\eco
\begin{proof} 
    Since $H$ is a basic vector field and 
     $h\na_UH=0$ for every $U$ vertical vector field,
     it follows $A_H\equiv 0$, by lemma \ref{lemdoi}.
     By dim \ho $=2$ and $H_p\not=0$ for every
     $p\in M$ we have $A\equiv 0$.
     Therefore, by proposition 9.104 in \cite{bes},
      $M$ is a locally warped product.
\end{proof}

We would like to know how much a semi-Riemannian submersion with totally
umbilic fibres is different to be a locally warped product. For  this
purpose we assume that $R(X,Y,X,Y)$ is constant along fibres 
for every $X$,$Y$ basic vector fields 
and the mean curvature vector field $H$ is parallel in the
horizontal bundle along fibres.  

First, we give equivalent conditions to these assumptions.

\bpr\label{ci}
   Let \pis\ be a semi-Riemannian submersion. 
   Then the following conditions are equivalent:\\
   i) $R(X,Y,X,Y)$ is constant along fibres 
   for every $X$, $Y$, $Z$ basic vector fields;\\
   ii) the function $g(A_XY,A_XZ)$ is constant along fibres
   for every $X$, $Y$, $Z$ basic vector fields;\\
   iii) $hR(X,Y)Z$ is a basic vector field
   for every $X$, $Y$, $Z$ basic vector fields.
\epr
\begin{proof}
$i)\Rightarrow ii)$\ Let $X$, $Y$ be basic vector fields. 
    By lemma \ref{lemtrei}, we have 
       $$ R(X,Y,X,Y)=R'(\pi_*X,\pi_*Y,\pi_*X,\pi_*Y)\circ\pi-3g(A_XY,A_XY)$$
   If $ R(X,Y,X,Y)$ is constant along fibres then
   $g(A_XY,A_XY)$ is constant along fibres. By polarization, we get
   $g(A_XY,A_XZ)$ is constant along fibres.\\
$ii)\Rightarrow iii)$\ 
   If $g(A_XY,A_XZ)$ is constant along fibres for every $X$, $Y$, $Z$
   basic vector fields then
   $A_XA_XY$ is basic for every $X$, $Y$  basic vector fields, by lemma \ref{lem}.
   Therefore, by polarization, $A_XA_YZ+A_YA_XZ$ is basic for every $X$, $Y$, $Z$
   basic vector fields.
   By O'Neill's equations (see \cite{one}), we have
     $hR(X,Y)Z=R'(X,Y)Z-2A_ZA_XY+A_XA_YZ-A_YA_XZ$,
   where $R'(X,Y)Z$ is the horizontal lifting of $R'(\pi_*X,\pi_*Y)\pi_*Z$.
   Rewriting this formula we get
     $$hR(X,Y)Z=R'(X,Y)Z-(A_ZA_XY+A_XA_ZY)+(A_ZA_YX+A_YA_ZX).$$
   Hence $hR(X,Y)Z$ is a basic vector field 
   for every $X$, $Y$, $Z$ basic vector fields.
\end{proof}

Let $\rho^{\mathcal V}(E)=\sum\limits_i\ei R(E,e_i)e_i$ 
for every $E$ tangent vector field to $M$.

\bpr\label{equival}
   Let \pis\ be a semi-Riemannian submersion with totally umbilic fibres.
   If $r\geq 2$ then the following conditions are equivalent:\\
   i) the mean curvature vector field $H$ is parallel in the horizontal bundle
     along fibres,\\
   ii) $\rho^{\mathcal V}(U)$ is vertical for every $U$ vertical vector field.
\epr
\begin{proof}
       By O'Neill's equation (see \cite{one}), we have
        $$R(e_i,U,e_i,X)=g((\na_UT)_{e_i}e_i,X)-g((\na_{e_i}T)_Ue_i,X).$$
      Then 
      \bea
           \sum\limits_i\ei g((\na_UT)_{e_i}e_i,X)
           &=&
                     \sum\limits_i\ei [g(\na_UT_{e_i}e_i,X)-
                                      g(T_{\na_U{e_i}}e_i,X)\no\\
                                      &&-g(T_{e_i}\na_U{e_i},X)]\no\\
          &=&
             g(\na_UH,X)-\sum\limits_i 2\ei g(\na_Ue_i,e_i)g(\f Hr,X)\no\\
          &=&
             g(\na_UH,X)-\sum\limits_i\ei U(g(e_i,e_i))\ g(\f Hr,X)\no\\
         &=&
             g(\na_UH,X)\no
        \eea
    Since the fibres are totally umbilic, we also have 
	$$\sum\limits_i\ei g((\na_{e_i}T)_{e_i}U,X)=\f 1rg(\na_UH,X).$$
    Replacing these in O'Neill's equation, we obtain
        $$g(\rho^{\mathcal V}(U),X)= \sum\limits_i\ei R(e_i,U,e_i,X)=
                (1-\f 1r)g(\na_UH,X).$$
    Therefore  $\rho^{\mathcal V}(U)$ is a vertical vector field 
    for every $U$ vertical vector field if and only if $H$ is parallel in the  
    horizontal bundle along fibres.
\end{proof}
\bth\label{tts}
    Let \pis\ be a semi-Riemannian submersion with totally umbilic fibres.
    We suppose that $R(X,Y,X,Y)$ is constant along fibres,
    $H$ is parallel in the horizontal bundle along fibres, 
    $s-s'\in\{0,r\}$. Then $H$ is a basic vector field if and only if
    $hR(X,Y)A_XY$ is a basic vector field for every $X$, $Y$ basic vector fields.
\eth
\begin{proof}
   {\it a})\ We suppose that $hR(X,Y)A_XY$ is a basic vector field 
             for every $X$, $Y$ basic vector fields.
             Let $X$, $Y$, $Z$ be basic vector fields.\\
    Using O'Neill's equation (see \cite{one}) we get
     \bea 
       R(X,Y,X,A_XY)&=&g((\na_XA)_XY,A_XY)+2g(A_XY,T_{A_XY}X)\no\\
                     &=&g(\na_X A_XY,A_XY))-g(A_{\na_XX}Y,A_XY)
                        -g(A_X\na_XY,A_XY)\no\\
                       &&+2g(A_XY,T_{A_XY}X)\no\\
                     &=&\f 12 X(g(A_XY,A_XY))-g(A_Y{h\na_XX},A_YX)
                        -g(A_Xh\na_XY,A_XY)\no\\
                       &&+2g(A_XY,T_{A_XY}X).\no
    \eea
    By proposition \ref{ci},
    the function $g(A_XY,A_XZ)$ is constant along fibres.
    If $hR(X,Y)A_XY$ is basic then $g(A_XY,T_{A_XY}X)$
    is constant along fibres.
    Let $U$ be a vertical vector field. Then
     \bea
           0=U(g(T_{A_XY}X,A_XY))= &-&\f 1r U(g(X,H))\ g(A_XY,A_XY)\no\\
                                    &-&\f 1r g(X,H)\ U(g(A_XY,A_XY)),\no
     \eea
      thus, we get
          \begin{equation}\label{xx} 0=U(g(X,H))\ g(A_XY,A_XY)\end{equation}
      for every $X$, $Y$ basic vector fields.\\
      Let $p$ in $M$ be an arbitrary point.
      Let $Y$ be a basic vector field such that $Y_p=H_p$.\\
      The relation (\ref{xx}) in point $p$ become
      \begin{equation}\label{xxx} 
          0=U_p(g(X,H))\ g_p(A_{X_p}H_p,A_{X_p}H_p)
      \end{equation}
      We have two possible situations:
\begin{case}
       $ A_{X_p}H_p\neq 0.$\\ 
       Since the metrics of fibres are negative definite for
       $s-s'=r$ or positive definite for
       $s-s'=0$ we have
       $g_p(A_{X_p}H_p,A_{X_p}H_p)\neq 0$.
       By relation (\ref{xxx}), we get
       $U_p(g(X,H))=0$.
\end{case}
\begin{case}
     $ A_{X_p}H_p=0$\\
     Since $H$ is parallel in the horizontal bundle
     along fibres we have $ U(g(H,X))=-g(A_XH,U)$.
     Using the hypothesis of {\it Case} 2,  $ A_{X_p}H_p=0$, we get
     $ U_p(g(H,X))=0$
\end{case}

     In both cases we proved that  $ U_pg(H,X)=0$
     for an arbitrary $p\in M$ and for every $X$ basic vector field.
     By lemma \ref{lem}, it follows that $H$ is a basic vector field.   

\vspace*{7pt}
{\it b})\ 
     We suppose $H$ is a basic vector field and we shall prove
     $hR(X,Y)A_XY$ is basic for every $X$, $Y$ basic vector fields.\\
     By O'Neill's equation, we get
\bea
       R(X,Y,Z,A_XY)&=&\f 12 Z(g(A_XY,A_XY))
                        -g(A_Y{h\na_ZX},A_YX)\no\\
                     && -g(A_Xh\na_ZY,A_XY)
                        -g(A_XY,A_XY)g(\f Hr,Z)\no\\
                     && -g(A_YZ,A_YX)g(\f Hr,X)
			-g(A_XZ,A_XY)g(\f Hr,Y).\no
\eea    
      Hence $g(hR(X,Y)A_XY,Z)$ is constant along fibres 
      for every $X$, $Y$, $Z$ basic vector fields. 
      Therefore, by lemma \ref{lem}, $hR(X,Y)A_XY$ is basic.
\end{proof}
      We apply the theorem  \ref{tts} in two particular cases,
      when the total space is either a constant 
      curvature semi-Riemannian manifold
      or a generalized complex
      space form.
\bco\label{tp}
    Let \pis\ be a semi-Riemannian submersion with totally umbilic fibres.
    If $M$ is a semi-Riemannian manifold with constant curvature and if 
    $r\geq 2$, $s-s'\in\{0,r\}$ then $H$ is a basic vector field and 
	$A_H\equiv 0$.
\eco
\begin{proof}
    If $M$ has constant sectional curvature $c$ then 
    for every $X$, $Y$ basic vector fields and for every $U$ vertical vector field
    we get:\\
    1) $R(X,Y,X,Y)=c(g(X,X)g(Y,Y)-g(X,Y)^2)$ is constant along fibres;\\
    2) $R(X,Y)A_XY=c(g(A_XY,X)X-g(A_XY,Y)Y)=0$;\\
    3)  $\rho^{\mathcal V}(U)=\sum\limits_i\ei R(U,e_i)e_i=
        \sum\limits_i\ei c(g(e_i,e_i)U-g(U,e_i)e_i)=
        c(r-1)U$ is a vertical vector field.
    Therefore, by theorem \ref{tts} and 
    proposition \ref{equival}, $H$ is a basic vector field.
    Hence, by lemma \ref{lemdoi}, $A_H\equiv 0$.
\end{proof}
    The totally umbilic submanifolds of 
    generalized complex space forms was classified by L. Vanhecke
    in \cite{van} (see also survey paper \cite{ianu}). Applying theorem \ref{tts}
    we get
\bco
    Let $\pi:(M,g)\to(B,g')$ be a semi-Riemannian 
    submersion with totally umbilic fibres
    from a generalized complex space form 
    onto an almost hermitian manifold. 
    If $\pi$ is a holomorphic map, $r\geq 2$, $s-s'\in\{0,r\}$ 
    then $H$ is a basic vector field and $A_H\equiv 0$.
\eco
\begin{proof}
    If $(M,J)$ is a generalized complex space form 
    of constant holomorphic sectional curvature
    $\mu$ and of type $\al$ then 
    the curvature tensor field satisfies (see \cite{ianu})
    \bea
    R(E,F,G,G')&=&\f 14(\mu+3\al)\{g(E,G)g(F,G')-g(E,G')g(F,G)\}\no\\
                &&+\f 14(\mu-\al)\{g(E,JG)g(F,JG')-g(E,JG')g(F,JG)\}\no\\
                &&+\f 12(\mu-\al)g(E,JF)g(G,JG').\no
    \eea
    Since $\pi$ is a holomorphic map, 
    we get $JX$ is basic for every $X$ basic vector field.
    Hence $R(X,Y,X,Y)$ is constant along fibres, $g(\rho^{\mathcal V}(U),X)=0$
    and $R(X,Y,Z,U)=0$ for every $X$, $Y$, $Z$ basic vector fields 
    and for every $U$ vertical vector field.
    Therefore, by theorem \ref{tts} and 
    proposition \ref{equival}, $H$ is a basic vector field.
\end{proof}
\bpr\label{ts}
    Let \pis\ be a semi-Riemannian submersion with totally umbilic fibres.
    If $R(X,Y,X,Y)$ is constant along fibres
    for every $X$, $Y$ basic vector fields,
    $H$ is basic and $B$ 
    is a compact and orientable manifold
    then $\ta$ is constant along fibres and
  \begin{equation}\label{sxxxx}
    \int\limits_B \tau'{ ^{HV}} =-\f 1r\int\limits_B g'(\pi_*H,\pi_*H)\mathrm{dv_g'}
                  +\int\limits_B g'(A,A)\mathrm{dv_g'}
  \end{equation}
    where $g'(A,A)$, $\tau'{ ^{HV}}$ are the functions on $B$ 
    satisfying $g'(A,A)\circ\pi=g(A,A)$ and  $\tau'{ ^{HV}}\circ\pi=\ta$.
\epr
\begin{proof}
     By theorem \ref{forjan} we have
      $$    \ta =\mathrm{div}(H)+(1-\f 1r)g(H,H)+g(A,A).$$
     Using lemma \ref{lemdoi} we get
 \bea
       \mathrm{div}(H)+g(H,H)&=&\sum\limits_\al\ea g(\na_{e_\al}H,{e_\al})+
       \sum\limits_i\ei g(\na_{e_i}H,{e_i})+g(H,H)\no\\
        &=&
       \sum\limits_\al\ea  g'(\na '_{\pi_*e_\al}\pi_*H,\pi_*{e_\al})\circ\pi+
       \sum\limits_i\ei [g(T_{e_i}H,e_i)+g(H,T_{e_i}e_i)]\no\\
       &=&\mathrm{div}'(\pi_*H)\circ\pi\no
  \eea
     Since the function $g(A_{e_\al}e_\bet,A_{e_\al}e_\bet)$ 
     is constant along fibres 
     we can consider the function $g'(A,A)$ on $B$ 
     given by $g'(A,A)\circ\pi=g(A,A)$.
     Then
  \begin{equation}\label{saptezece}
     \ta=[\mathrm{div}'(\pi_*H)-\f 1r g'(\pi_*H,\pi_*H)+g'(A,A)]\circ\pi.
  \end{equation}
  and $\ta$ is constant along fibres. Let $\tau'{ ^{HV}}$ be the function on 
  $B$ such that $\tau'{ ^{HV}}\circ\pi=\ta$.
     Since $\int\limits_B \mathrm{div}'(\pi_*H)=0$, it follows the formula
             (\ref{sxxxx}).
\end{proof}
\bco\label{treisase}
    Let \pis\ be a semi-Riemannian submersion with totally umbilic fibres.
    If $M$ is a semi-Riemannian manifold with constant curvature $c$,
    $r\geq 2$, $s-s'\in\{0,r\}$  and $B$ is a compact and orientable manifold  then
    \begin{equation}\label{xxxx}
    rnc\mathrm{vol}(B)=-\f 1r\int\limits_B g'(\pi_*H,\pi_*H)\mathrm{dv_g'}
                  +\int\limits_B g'(A,A)\mathrm{dv_g'}
    \end{equation}
    where $g'(A,A)$ is the function on $B$ satisfying $g'(A,A)\circ\pi=g(A,A)$.
\eco
\begin{proof}
     Since $M$ is a semi-Riemannian manifold with constant curvature $c$,
     we get $R(e_\al,e_i,e_\al,e_i)=\ei\ea c$. So $\ta=rnc$.
\end{proof}
\bth\label{t:t}
     Let $B$ be an $n$-dimensional compact and orientable Riemannian manifold,
     $\pi :M\to B$ be a semi-Riemannian submersion with totally umbilic fibres.
    We suppose that $M$ is an $m$-dimensional  semi-Riemannian manifold of index $r=m-n$
    with non-negative mixed curvature, $R(X,Y,X,Y)$ is constant along fibres,
    for every $X$, $Y$ basic vector fields, and
    $H$ is basic.
    Then\\
    i)\ $K(X,V)=0$ for every $X$ horizontal vector field 
        and for every  $V$ vertical vector field;\\
    ii)\ the fibres are totally geodesic and the horizontal distribution is integrable,
        hence $M$ is a locally warped product.
\eth
\begin{proof}
     Since $\ea =1$, the metrics induced on fibres are negative definite, we have 
     $g(A,A)=\sum\limits_{\al\ ,\ \bet}
           \ea\eb g(A_{e_\al}e_\bet , A_{e_\al}e_\bet ) \leq 0$.
     By proposition \ref{ts}, $\int\limits_B\tau'{^{ HV}}\leq 0$.\\
     But the mixed curvature is non-negative, hence $\tau'{^{ HV}}\geq 0$.
     Therefore $K(X,V)=0$  for every $X$ horizontal vector field 
        and for every  $V$ vertical vector field and
     $A\equiv 0$, $T\equiv 0$.
\end{proof} 
\bco\label{treisapte}
     Let $B$ be an $n$-dimensional compact and orientable Riemannian manifold,
     $\pi :M\to B$ be a semi-Riemannian submersion with totally umbilic fibres.
    We suppose $M$ is an $m$-dimensional  semi-Riemannian manifold of index $r=m-n$
    with constant curvature $c$ and $r\geq 2$.
    Then\\
    i)\ $c\leq 0$;\\
    ii)\ If $c=0$ then the fibres are totally geodesic and the
    horizontal distribution \ho\ is integrable.
\eco
\begin{proof}
     If we suppose $c>0$ then, by theorem \ref{t:t}, we get $c=K(X,V)=0$,
     which is a contradiction with our assumption. Therefore $c\leq 0$.
\end{proof} 
\bco
      Let $B$ be an $n$-dimensional simply connected Riemannian manifold,
     $\pi :M\to B$ be a semi-Riemannian submersion with totally umbilic fibres.
     If $M$ is an $m$-dimensional  semi-Riemannian manifold of index 
     $r=m-n\geq 2$ with constant curvature $c$ then $B$ is not compact.
\eco
\begin{proof}
     If we suppose that $B$ is a compact Riemannian manifold then $B$ is complete.
     By corollary \ref{treisapte}, we get $c\leq 0$. 
     It follows $K'\leq 0$, by lemma \ref{lemtrei}.
     By  Hadamard's theorem, we have $B$ is diffeomorphic to $\mathbb R^n$.
     Therefore $B$ is not compact.
\end{proof}     
     
We introduce the notion of Clairaut semi-Riemannian submersion (see \cite{all}).

\begin{defn}
     Let $B$ be an $n$-dimensional Riemannian manifold,
     $\pi :(M,g)\to (B,g')$ be a semi-Riemannian submersion,
     M be an $m$-dimensional semi-Riemannian 
     manifold of index $r=m-n$ and $\ga $ be a timelike geodesic in $M$.
    We denote the velocity vector field of $\ga$ by $E=\ga'$ and its vertical
    part by $V$.
    At each point $\ga(s)$ we define $\varphi(s)$ to be the hyperbolic angle
    between $E$ and $V$, i.e. $\varphi\geq 0$ is the number satisfying:
     $$g(E,V)=-\|E\|\cdot \|V\| \cosh\varphi ,$$
    where $\|E\|^2=-g(E,E)$ and $\|V\|^2=-g(V,V)$.\\
    $\pi :(M,g)\to (B,g')$
    is said to be a {\it Clairaut semi-Riemannian submersion} 
   if there is a positive
   function $w:M\to \Bbb R$ such that for every timelike geodesic $\ga$ in $M$,
   $w\cosh\varphi$ is constant along $\ga$. We call $r$ the {\it girth} of
    the submersion.
\end{defn}

	The following characterization of Clairaut semi-Riemannian submersion
   given by D.Allison \cite{all} for $1$-dimensional fibres has a similar proof for
   $r$-dimensional case
\bpr\label{treiopt}
    Let $B$ be an $n$-dimensional Riemannian manifold,
    $\pi :(M,g)\to (B,g')$ be a semi-Riemannian submersion
    and $M$ be an $m$-dimensional semi-Riemannian 
    manifold of index $r=m-n$.
   Then $\pi$ is a Clairaut semi-Riemannian submersion
   with girth $w=\mathrm{exp } f$ if and only if the fibres are totally umbilic
   with a gradient $H=-\mathrm{grad }f$ as mean curvature vector field.
   Furthermore for a Clairaut  semi-Riemannian submersion
   having connected fibres $f=f_*\circ\pi$ for some $f_*:B\to\Bbb R$
   and $H$ is a basic vector field obtained by lifting $H_*=-\mathrm{grad }f_*$
   horizontally. 
\epr
\ble\label{treinoua}
    Let $\pi :(M,g)\to (B,g')$ be a semi-Riemannian submersion 
    with totally umbilic fibres. We suppose 
    that $H$ is a basic vector field and 
    parallel in the horizontal bundle along fibres
    and there is a constant $c\in\Bbb R$ such that
    $g(\rho^{\mathcal V}(\f Hr),X)=cg(H,X)$ 
    for every $X$ horizontal vector field.
    Then:
\begin{equation}\label{lx}
   \f 12\mathrm{grad}(c+g(\f Hr,\f Hr))=(c+g(\f Hr,\f Hr))\f Hr
\end{equation}
\ele
\begin{proof}
  Let $X$,$Y$ be horizontal vector fields and $U$ be a vertical vector field.
  By O'Neill's equation \cite{one} we have
  \bea
  R(X,U,Y,U)&=&g(U,U)[g(\na_X\f Hr,Y)-g(X,\f Hr)g(Y,\f Hr)]\no\\
            &&+g((\na_UA)_XY,U)+g(A_XU,A_YU)\no
  \eea
  Let $Y=\f Hr$. Since $H$ is a basic vector field and 
    parallel in the horizontal bundle along fibres we get 
    $A_H\equiv 0$. Therefore
  $$g((\na_UA)_XH,U)+g(A_XU,A_HU)=0$$
  Since 
       $$\sum\ei R(X,e_i,\f Hr,e_i)=g(\rho^{\mathcal V}(\f Hr),X)=cg(H,X)$$
  we obtain
  \bea
 &&  g(X,\f Hr)c=g(\na_X\f Hr,\f Hr)-g(X,\f Hr)g(\f Hr,\f Hr);\no\\
 &&
  g(X,\f Hr)(c+g(\f Hr,\f Hr))=g(\na_X\f Hr,\f Hr)=
    \f 12 Xg(\f Hr,\f Hr);\no\\
 &&
  g(X,\f 12\mathrm{grad}\ g(\f Hr,\f Hr)-
  (c+g(\f Hr,\f Hr))\f Hr)=0\no
\eea
   for every $X$ horizontal vector field.\\
   This implies
   $\f 12h\ (\mathrm{grad}\ g(\f Hr,\f Hr))=
   (c+g(\f Hr,\f Hr))\f Hr$.\\
   Since $H$ is parallel in the horizontal bundle along fibres, it follows that
   $g(H,H)$ is constant along fibres. Hence
   $v\mathrm{grad}(g(\f {H}m,\f {H}m)+c)=0$.
\end{proof}
\ble\label{lemnew}
   Let $\pi:M\to B$ be a semi-Riemannian submersion. If either\\
   \br i)\er\ $M$ has constant curvature, \ or\\
   \br ii)\er\ $M$ and $B$ are Einstein manifolds and $A_H\equiv 0$\\
   then there is  $c\in\Bbb R$ such that 
    $$g(\rho^{\mathcal V}(\f Hr),X)=cg(H,X)$$
   for every $X$ horizontal vector field. 
\ele
\begin{proof}
  i)  If $M$ has constant curvature $c$ then
     $$g(\rho^{\mathcal V}(\f Hr),X)=\sum\limits_i\ei R(\f Hr,e_i,X,e_i)=cg(H,X)$$
  for every $X$ horizontal vector field. \\
  ii) We define $$\rho^{\mathcal H}(E)=\sum\limits_\al\ea R(E,e_\al)e_\al$$
      for every $E$ vector field on $M$. We have
      $$\rho(\f Hr,X)=g(\rho^{\mathcal V}(\f Hr),X)+g(\rho^{\mathcal H}(\f Hr),X).$$
      Since $M$ and $B$ are Einstein manifolds, for some $\la$, $\la'\in\mathbb R$ we have
      $\rho(\f Hr,X)=\la g(\f Hr,X)$ and 
       $\rho'(\pi_*\f Hr,\pi_*X)=\la' g(\pi_*\f Hr,\pi_*X)$
      for every $X$ horizontal vector.\\
     By lemma \ref{lemtrei}, we get
      \bea
          \rho'(\pi_*\f Hr,\pi_*X) &=&
              \sum\limits_\al\ea(R(\f Hr,e_\al,X,e_\al)+3g(A_{e_\al}\f
						Hr,A_{e_\al}X))\no\\
       				   &=& g(\rho^{\mathcal H}(\f Hr),X)\no
     \eea
     It follows
       $$g(\rho^{\mathcal V}(\f Hr),X)=(\la-\la')g(\f Hr,X)$$
   for every $X$ horizontal vector field.
\end{proof}
\bth\label{nxzece}
    Let
    $\pi :(M,g)\to (B,g')$ be a semi-Riemannian submersion 
    with totally umbilic fibres. We suppose
    that $H$ is a basic vector field and 
    parallel in the horizontal bundle along fibres
    and there is a constant $c\in\Bbb R$ such that
    $g(\rho^{\mathcal V}(\f Hr),X)=cg(H,X)$ for every $X$ horizontal vector field.
    If $g(\f {H_p}r,\f {H_p}r)\neq -c$ for all $p\in M$
    then
    $$\f Hr=\f 12\mathrm{grad(ln}|c+g(\f Hr,\f Hr)|).$$
    If moreover $B$ is an $n$-dimensional Riemannian manifold,
     and M be an $m$-dimensional semi-Riemannian 
    manifold of index $r=m-n$ then
    $\pi$ is a Clairaut semi-Riemannian submersion.
\eth
\begin{proof}
   Since $g(\f {H_p}r,\f {H_p}r)+c\neq 0$ for every point $p\in M$,
   we have, by formula \eqref{lx},
   $$\f Hr=\f 12\mathrm{grad(ln\ }|g(\f Hr,\f Hr)+c|).$$
   By proposition \ref{treiopt}, we get $\pi$ is a 
   Clairaut semi-Riemannian submersion.
\end{proof}
\bco\label{xzece}
    Let $B$ be an $n$-dimensional Riemannian manifold,
    $\pi :(M,g)\to (B,g')$ be a semi-Riemannian submersion 
    with totally umbilic fibres
    and M be an $m$-dimensional semi-Riemannian 
    manifold of index $r=m-n\geq 2$ with constant curvature $c$.
    If $g(\f {H_p}r,\f {H_p}r)\neq -c$ for all $p\in M$
    then
    $$\f Hr=\f 12\mathrm{grad(ln}|c+g(\f Hr,\f Hr)|)$$
    and $\pi$ is a Clairaut semi-Riemannian submersion.
\eco
\begin{proof}
    By lemma \ref{lemnew} and by theorem \ref{nxzece}
    we get the conclusion.
\end{proof}

    By corollaries \ref{tp} and  \ref{treisase} and use of \cite{kim},
    we get more information about
    curvature of total space,
    in the constant curvature case.
\bpr    
        Let $\pi :(M,g)\to (B,g')$ be a semi-Riemannian submersion 
    with totally umbilic fibres. We assume that
    $B$ is an $n$-dimensional Riemannian manifold,
    $M$ is an $m$-dimensional semi-Riemannian 
    manifold of index\\ $r=m-n\geq 2$ with constant curvature $c$.
    If $g(\f Hr,\f Hr)+c\geq 0$ everywhere,
    $g(\f {H_{p_0}}r,\f {H_{p_0}}r)+c>0$ at some point $p_0\in M$,
    and if each fibre is a compact manifold then
    the horizontal distribution is integrable and $M$ is a locally warped
    product.
    If moreover $B$ is a compact and orientable manifold then
    $c<0$, $n\neq 1$ and $M$ is not a compact and orientable manifold.
\epr
\begin{proof}
    If $M$ has constant curvature, then for every $U$ vertical vector field
    we have
    $\hat\rho(U,U)=(r-1)g(U,U)(g(\f Hr,\f Hr)+c)$.
    Since the metrics of the fibres are negative definite we have
    $\hat\rho(U,U)\leq 0$ everywhere and $\hat\rho$ is negative definite in $p_0$.    
    By corollary \ref{tp}, we have $H$ is a basic vector field. Therefore,
    by the theorem in \cite{kim}, the horizontal distribution is integrable
    and $M$ is a locally warped product.
    If $B$ is compact and orientable then, by corollary \ref{treisase},
    $$(n-1)c\mathrm{vol}(B)=
        -\int\limits_B(g'(\f {\pi_*H}r,\f {\pi_*H}r)+c))\mathrm{dv_{g'}}.$$
    Hence $(n-1)c<0$.
    If we suppose $M$ is a compact and orientable manifold 
    then by proposition \ref{for} we have
    $$(n+r-1)c\mathrm{vol}(M)=(r-1)\int\limits_M (g(\f Hr,\f Hr)+c)
                                                     \mathrm{dv_{g}}
                          +\f 1r\int\limits_M g(A,A)\mathrm{dv_{g}}$$
    By $g(\f Hr,\f Hr)+c>0$ in $p_0$ and $g(A,A)\equiv 0$ we get $c>0$, 
    which is  a contradiction
    with  corollary \ref{treisapte}.
\end{proof}
\bpr\label{unspe}
    Let $B$ be an $n$-dimensional compact and orientable Riemannian manifold,
    $\pi :(M,g)\to (B,g')$ be a semi-Riemannian submersion 
    with totally umbilic fibres
    and M be an $m$-dimensional semi-Riemannian 
    manifold of index $r=m-n$ with constant curvature $c$.
    If  each fibre is a compact manifold, $r\geq 2$, $n\geq 3$, and if 
        \begin{equation}\label{l}
           g(\f {H_p}r,\f {H_p}r)+c\geq 0
        \end{equation}
        \begin{equation}\label{ll}
           3 g(\f {H_p}r,\f {H_p}r)+nc-\f 1r g_p(A,A)\geq 0
        \end{equation}
          for all $p\in M$ then  $g(\f Hr,\f Hr)+c\equiv 0$.
\epr
\begin{proof}
     We denote by 
         $\De_B f=\mathrm{div(grad\ }f)$ the Laplacian of $B$ defined on
    functions $f:B\to \Bbb R$.

    Let $f$ be the function on $B$ given by
        $$f=g'(\f {\pi_*H}r,\f {\pi_*H}r)+c$$
    By lemma \ref{treinoua} we have the equation 
            $\f 12 \mathrm{grad\ }f=f \f {\pi_*H}r$.
    We get
\bea
     \f 12 \De_B f&=&\f 12\mathrm{div(grad\ }f)=\mathrm{div}(f \f {\pi_*H}r)=
      (\f {\pi_*H}r)(f)+f\mathrm{div}(\f {\pi_*H}r)\no\\
         &=& g'(\mathrm{grad} f,\f {\pi_*H}r)+f\mathrm{div}(\f {\pi_*H}r)\no\\
        &=&2fg'(\f {\pi_*H}r,\f {\pi_*H}r)+f\mathrm{div}(\f {\pi_*H}r)\no
\eea
     Using relation (\ref{saptezece}) we have 
    $$\mathrm{div}(\f {\pi_*H}r)=nc+g'(\f {\pi_*H}r,\f {\pi_*H}r)
           -\f 1rg'(A,A),$$
    where $g'(A,A)$ is the function satisfying $g(A,A)=g'(A,A)\circ\pi$.
    We obtain 
\bea
    \f 12 \De_B f&=& f(3(f-c)+nc-\f 1r g'(A,A))\no\\
                 &=& (g'(\f {\pi_*H}r,\f {\pi_*H}r)+c)
                      (3g'(\f {\pi_*H}r,\f {\pi_*H}r)+nc-\f 1rg'(A,A))\geq 0\no
\eea
    Since $\De_B (f)\geq 0$ and $B$ is a compact and orientable manifold,
     we have $f$ is a constant function and  $\De_B (f)\equiv 0$, by Hopf's lemma.

    If we suppose $f=g'(\f {\pi_*H}r,\f {\pi_*H}r)+c>0$ then, by corollary
    \ref{xzece}, $A\equiv 0$, $c<0$.\\ 
    $\De_B f\equiv 0$ implies
    $3g'(\f {\pi_*H}r,\f {\pi_*H}r)+nc\equiv 0$.
    From corollary \ref{treisase}, we obtain
       $\int\limits_B (g'(\f {\pi_*H}r,\f {\pi_*H}r)+nc)\mathrm{dv}_{g'}=0$. Since
   $g'(\f {\pi_*H}r,\f {\pi_*H}r)$ is a nonnegative constant function we have
   $g'(\f {\pi_*H}r,\f {\pi_*H}r)=0$. It follows $H\equiv 0$. Therefore
   $c>0$. This is a contradiction with condition $c\leq 0$ given by corollary
  \ref{treisapte}.

  All these imply $g'(\f {\pi_*H}r,\f {\pi_*H}r)+c\equiv 0$.
\end{proof}
\ble\label{lemz}
   The condition $ 3 g(\f {H_p}r,\f {H_p}r)+nc\geq 0$
   for every $p\in M$ implies the conditions 
   \eqref{l} and \eqref{ll}
   in proposition \br \ref{unspe}\er .
\ele
\begin{proof} Since $g(A,A)\leq 0$ and 
        $ (n-3)g(\f {H_p}r,\f {H_p}r)\geq 0$
   for every $p\in M$ we have the conditions
  \eqref{l} and \eqref{ll}
   in proposition \ref{unspe}.
\end{proof}
By lemma \ref{lemz} and proposition \ref{unspe}, we have the following corollary.
\bco
    Let $B$ be an $3$-dimensional compact and orientable Riemannian manifold,
    $\pi :(M,g)\to (B,g')$ be a semi-Riemannian submersion 
    with totally umbilic fibres
    and M be an $m$-dimensional semi-Riemannian 
    manifold of index $r=m-3$ with constant curvature $c$.
    If  $r\geq 2$, if each fibre
    is a compact manifold and if
    $g(\f Hr,\f Hr)+c\geq 0$  then 
    $g(\f Hr,\f Hr)+c\equiv 0$.
\eco
\begin{acknow}
The second author wishes to thank Professors M. Falcitelli and
A.M. Pastore for useful discussions on Riemannian submersions
during his stay in the Autumn of 1999 
as Visiting Professor in Bari. The second author benefits of financial support from 
    EDGE-project.
\end{acknow}

\end{document}